\documentclass[11pt]{amsart}

\usepackage{graphicx,calc}

\newtheorem{theorem}{Theorem} 

\newtheorem{lemma}[theorem]{Lemma}

\newcommand{\longpage}{\enlargethispage{\baselineskip}}

\newcommand{\mapright}[1]{\stackrel{#1}{\longrightarrow}}

\begin{document}

\title{Imbeddings of free actions on handlebodies}

\author{Darryl McCullough}
\address{Department of Mathematics\\
University of Oklahoma\\
Norman, Oklahoma 73019\\
USA}
\email{dmccullough@math.ou.edu}
\urladdr{www.math.ou.edu/$_{\widetilde{\phantom{n}}}$dmccullo/}
\thanks{The author was supported in part by NSF grant DMS-0102463}
\subjclass{Primary 57M50; Secondary 57M60, 58D99}

\subjclass{Primary 57M60; Secondary 20F05}

\date{\today}

\keywords{3-manifold, handlebody, group action, free, free action, imbed,
imbedding, equivariant, invariant, hyperbolic, Seifert, Heegaard, Heegaard
splitting, Whitehead link}

\begin{abstract}
Fix a free, orientation-preserving action of a finite group $G$ on a
$3$-dimensional handlebody $V$. Whenever $G$ acts freely preserving
orientation on a connected $3$-manifold $X$, there is a $G$-equivariant
imbedding of $V$ into $X$. There are choices of $X$ closed and
Seifert-fibered for which the image of $V$ is a handlebody of a Heegaard
splitting of $X$. Provided that the genus of $V$ is at least $2$, there are
similar choices with $X$ closed and hyperbolic.
\end{abstract}

\maketitle

\section*{Introduction}
\label{sec:intro}

Any finite group acts (smoothly and) freely preserving orientation on some
$3$-dimensional handlebody, and the number of inequivalent actions of a
fixed $G$ on a fixed genus of handlebody can be arbitrarily large
\cite{MW}. Consequently, the following imbedding property of such actions
may appear surprising at first glance:

\begin{theorem} Let $G$ be a finite group acting freely and preserving
orientation on two handlebodies $V_1$ and $V_2$, not necessarily of the
same genus. Then there is a $G$-equivariant imbedding of $V_1$ into $V_2$.
\end{theorem}

\noindent In fact, this result is almost a triviality, as is the following
theorem of which it is a special case:

\begin{theorem} Let $G$ be a finite group acting freely and preserving
orientation on a handlebody $V$ and on a connected $3$-manifold $X$. Then
there is a $G$-equivariant imbedding of $V$ into $X$.
\label{thm:always}
\end{theorem}

By a result of D. Cooper and D. D. Long \cite{CL}, any finite group acts
freely on some hyperbolic rational homology $3$-sphere. So
theorem~\ref{thm:always} shows that a free $G$-action on a handlebody
always has an extension to an action on such a $3$-manifold. Also, by a
result of S. Kojima \cite{K}, for any finite $G$ there is a closed
hyperbolic $3$-manifold whose full isometry group is $G$, and Kojima's
construction actually produces a free action. So there is an extension to a
free action on a closed hyperbolic $3$-manifold whose full isometry group
is~$G$.

One might ask for a more natural kind of extension, to a free $G$-action on
a closed $3$-manifold $M$, for which $V$ is one of the handlebodies in a
$G$-invariant Heegaard splitting of $M$. Simply by forming the double of
$V$ and taking an identical action on the second copy of $V$, one obtains
such an extension with $M$ a connected sum of $S^2\times S^1$'s. A better
question is whether $V$ is an invariant Heegaard handlebody for a free
action on an irreducible $3$-manifold. Our main result answers this
affirmatively.
\longpage

\begin{theorem} Let $G$ be a finite group acting freely and preserving
orientation on a handlebody $V$. Then the action is the restriction of a
free $G$-action on a closed irreducible $3$-man\-i\-fold $M$, which has a
$G$-invariant Heegaard splitting with $V$ as one of the handlebodies. One
may choose $M$ to be Seifert-fibered. Provided that $V$ has genus greater
than $1$, one may choose $M$ to be hyperbolic. In both cases, there are
infinitely many choices of $M$.\par
\label{thm:irreducible extension}
\end{theorem}
We remark that any orientation-preserving action of a finite group on a
closed $3$-manifold, free or not, has an invariant Heegaard splitting. For
the quotient is a closed orientable $3$-orbifold with $1$-dimensional
(possibly empty) singular set. One may triangulate the quotient so that the
singular set is a subcomplex of the $1$-skeleton. Then, the preimage of a
regular neighborhood of the $1$-skeleton is invariant and is one of the
handlebodies in a Heegaard splitting.

In the remaining sections of this paper, we prove theorems~\ref{thm:always}
and~\ref{thm:irreducible extension}. In \cite{MW}, a number of results
about free $G$-actions on handlebodies are obtained using more algebraic
methods.

\section{Proof of theorem~\ref{thm:always}}
\label{sec:always}

Recall that two $G$-actions on spaces $X$ and $Y$ are \emph{equivalent} if
there is a homeomorphism $j\colon X\to Y$ such that $h(x)=j^{-1}(h(j(x)))$
for all $x\in X$ and all $h\in G$. If $G$ acts properly discontinuously and
freely on a path-connected space $X$, then the quotient map $X\to X/G$ is a
regular covering map, so by the theory of covering spaces the action
determines an extension
\begin{equation*}
1\longrightarrow\pi_1(X)\longrightarrow\pi_1(X/G)
\mapright{\phi}G\longrightarrow1\ .
\end{equation*}
\noindent
Since we have not specified basepoints, the homomorphism $\phi$ is
well-defined only up to an inner automorphism of $G$.

Suppose now that $G$ is finite and acts freely and preserving orientation
on a handlebody $V$. The quotient manifold $V/G$ is orientable and
irreducible with nonempty boundary, so $\pi_1(V/G)$ is torsionfree. A
torsionfree finite extension of a finitely generated free group is free
(by~\cite{K-P-S} any finitely generated virtually free group is the
fundamental group of a graph of groups with finite vertex groups, and if
the group is torsionfree, the vertex groups must be trivial). So
$\pi_1(V/G)$ is free, and theorem~5.2 of~\cite{Hempel} shows that $V/G$ is
a handlebody. In this context, we obtain a simple algebraic criterion for
equivalence.
\begin{lemma} Suppose that $G$ acts freely and preserving orientation on 
handlebodies $V_1$ and $V_2$, with quotient handlebodies $W_1$ and $W_2$,
determining homomorphisms $\phi_i\colon \pi_1(W_i)\to G$. The actions are
equivalent if and only if there is an isomorphism $\Psi\colon \pi_1(W_1)\to
\pi_1(W_2)$ for which $\phi_2\circ \Psi=\phi_1$.
\label{lem:equivalence}
\end{lemma}
\begin{proof} An equivalence of the actions $j\colon V_1\to V_2$ induces a
homeomorphism $\overline{j}\colon W_1\to W_2$ for which $\phi_2\circ
\overline{j}_\# = \phi_1$. Conversely, suppose $\Psi$ exists. Since
both $W_1$ and $W_2$ are orientable, there is a homeomorphism $f\colon
W_1\to W_2$. Using well-known constructions of homeomorphisms of $W_2$ (as
for example in~\cite{M-M1}), all of Nielsen's \cite{Nielsen} generators of
the automorphism group of the free group $\pi_1(W_2)$ can be induced by
homeomorphisms, so $f$ may be selected to induce $\Psi$. The condition that
$\phi_2\circ \Psi=\phi_1$ then shows that $f$ lifts to a homeomorphism of
covering spaces $j\colon V_1\to V_2$, and moreover ensures that
$h(x)=j^{-1}(h(j(x)))$.
\end{proof}

Now we prove theorem~\ref{thm:always}. Let $W$ be the quotient handlebody
of the action on $V$, let $Y=X/G$, and let $\phi\colon \pi_1(W)\to G$ and
$\psi\colon \pi_1(Y)\to G$ be the homomorphisms determined by the actions.

There is an imbedding $k\colon W\to Y$ so that $\psi\circ k_\#=\phi$. For
we can regard $W$ as a regular neighborhood of a $1$-point union $K$ of
circles, so that $\pi_1(K)=\pi_1(W)$, and construct a map $k_0$ of $K$ into
$Y$ for which $\psi\circ (k_0)_\#=\phi$. Since $K$ is $1$-dimensional,
$k_0$ is homotopic to an imbedding, and since $W$ and $Y$ are orientable,
this imbedding extends to an imbedding $k$ of $W$ into $Y$. Since
$\psi\circ k_\#=\phi$, the preimage of $k(W)$ in $X$ is connected, and by
lemma~\ref{lem:equivalence} the restricted $G$-action on it is equivalent
to the original action on~$V$.

Theorem~\ref{thm:always} extends to the case when some elements of $G$
reverse the orientation. The equivariant imbedding exists if and only if
the subgroups of elements of $G$ that reverse orientation on $V$ and on $X$
are identical. The proof is affected only at the step when the imbedding of
$K$ into $Y$ is extended to an imbedding of $W$ into $Y$. The equality of
the orientation-reversing subgroups is precisely the condition needed for
the extension to exist.

\section{Proof of theorem~\ref{thm:irreducible extension}}

If the genus of $V$ is $0$, then $G$ is trivial and we take $M=S^3$. If the
genus of $V$ is $1$, so that $V=D^2\times S^1$, then (using
lemma~\ref{lem:equivalence}) every free $G$-action is equivalent to a
cyclic rotation in the $S^1$-factor. Regarding $V$ as a trivially fibered
solid torus in the Hopf fibering of $S^3$, the action extends to a free
action on $S^3$ with $V$ an invariant Heegaard splitting (it also extends
to free actions on infinitely many lens spaces containing $V$ as a fibered
Heegaard torus). So we may assume that the genus of $V$ is greater than
$1$. The quotient handlebody $W=V/G$ has genus at least~$2$ (since $V$ and
consequently $W$ have negative Euler characteristic).

\begin{figure}
\includegraphics[width=58 ex]{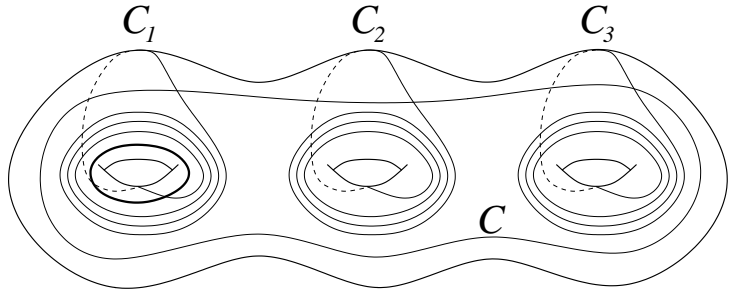}
\caption{The quotient handlebody $W$.}
\label{fig:Heegaard}
\end{figure}

We first construct the Seifert-fibered extension. As in section
\ref{sec:always}, there is a homomorphism $\phi\colon \pi_1(W)\to G$ that
determines the action. Let $g$ be the genus of $W$, and let $n$ be any
positive integer divisible by the orders of all the elements of $G$. We
consider a collection of simple closed curves $C_1,\dots\,$, $C_g$ in the
boundary $\partial W$, as shown in figure~\ref{fig:Heegaard} for the case
when $g=3$ and $n=4$. Each $C_i$ winds $n$ times around one of the handles
of $W$. Let $C_i'$ be the image of $C_i$ under the $n^{th}$ power of a Dehn
twist of $\partial W$ about the curve $C$. The union of the $C_i$ does not
separate $\partial W$, so neither does the union of the $C_i'$. So we can
obtain a closed $3$-manifold $Y$ with $W$ as a Heegaard handlebody by
attaching $2$-handles along the $C_i'$ and filling in the resulting
$2$-sphere boundary component with a $3$-ball.

Let $x_1,\ldots\,$, $x_g$ be a standard set of generators of $\pi_1(W)$,
where $x_i$ is represented by a loop that goes once around the $i^{th}$
handle. In $\pi_1(W)$, $C_i$ represents $x_i^n$ (up to conjugacy), and
$C_i'$ represents $x_i^n\,(x_1\cdots x_g)^{-n}$. Since every element of $G$
has order dividing $n$, it follows that $\phi$ carries each $C_i'$ to the
trivial element of $G$, so induces a homomorphism $\psi\colon \pi_1(Y)\to
G$. If $k\colon W\to Y$ is the inclusion, then $\psi\circ k_\#=\phi$. The
covering space $M$ of $Y$ has a free $G$-action and an invariant Heegaard
splitting, one of whose handlebodies is the covering space of $W$
corresponding to the kernel of $\phi$, that is, $V$.

\begin{figure}
\includegraphics[width=75 ex]{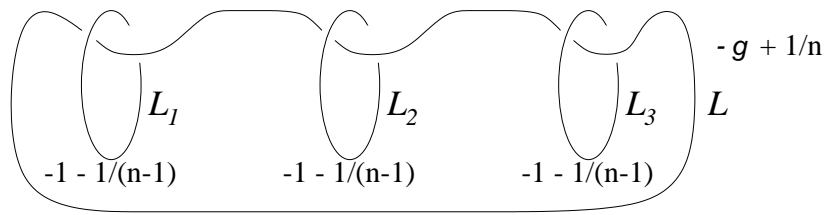}
\caption{Surgery description of $Y$ in the Seifert-fibered case.}
\label{fig:Dehn}
\end{figure}

We will show that $Y$ is Seifert-fibered, from which it follows that $M$ is
Seifert-fibered. Choose imbedded loops in the interior of $W$: $L$ near and
parallel to $C$, and $L_1,\dots\,$, $L_g$ near and parallel to loops
$\ell_i$ in $\partial W$ with each $\ell_i$ going once around the $i^{th}$
handle, meeting $C_i$ in one point. The loop $\ell_1$ appears in
figure~\ref{fig:Heegaard}. By a standard procedure, as explained for
example on pp.~275-278 of \cite{Rolfsen}, whose notation we follow, we may
change the attaching curves for the discs by Dehn twists about $C$ and the
$\ell_i$, at the expense of performing Dehn surgery on $L$ and the
$L_i$. First we twist $n$ times along $C$, introducing a $1/n$ coefficient
on $L$ and moving each $C_i'$ back to $C_i$. Then, $n-1$ twists along each
$\ell_i$ move $C_i$ to a loop $C_i''$ in $\partial W$ that looks like $C_i$
except it goes only once around the handle. This creates surgery
coefficients of $-1/(n-1)$ on the $L_i$. We may change the attaching
homeomorphism of the Heegaard splitting by any homeomorphism of $\partial
W$ that extends over $W$, without changing $Y$. In particular, we may
perform left-hand twists in the meridinal $2$-discs of the $2$-handles to
move the $C_i''$ to the $\ell_i$. This subtracts $1$ from the surgery
coefficients of the $L_i$, and subtracts $g$ from the coefficient of $L$,
yielding the diagram in figure~\ref{fig:Dehn}. The $\ell_i$ are the
attaching curves for the discs of a Heegaard description of $S^3$, so $Y$
is obtained from $S^3$ by Dehn surgery using the diagram.

We can already see a Seifert fibering, but we will work out the exact
Seifert invariants. The complement of a regular neighborhood of $L$ is a
solid torus $T$, for which the $L_i$ are fibers of a product fibering. A
cross-sectional surface in this fibering contained in a meridian disc of
$T$ meets the boundary of a regular neighborhood of each $L_i$ in a
meridian circle and meets the boundary of a regular neighborhood of $L$ in
the negative of the longitude, while the fiber meets them in longitude
circles and a meridian circle respectively.

A surgery coefficient $a/b$ means that a solid torus is filled in so that
$a\cdot m+b\cdot \ell$ becomes contractible, where $m$ and $\ell$ are a
meridian-longitude pair for a boundary torus of a regular neighborhood of
the link component. A Seifert invariant $(\alpha,\beta)$ determines a
filling in which $\alpha\cdot q+\beta\cdot t$ becomes contractible, where
$q$ is the cross-section and $t$ is the fiber. So the surgery coefficients
of our link produce one exceptional fiber with Seifert invariant $(n,1-n)$
for each $L_i$, and one with Seifert invariants $(n,gn-1)$ for $L$. In the
notation of \cite{Orlik}, the unnormalized Seifert invariants of $Y$ are
$\{0;(o_1,0);(n,1-n),\dots,(n,1-n),(n,gn-1)\}$, where there are $g+1$
exceptional orbits. The normalized invariants are
$\{-1;(o_1,0);(n,1),\dots,(n,1),(n,n-1)\}$.

\begin{figure}
\includegraphics[width=35 ex]{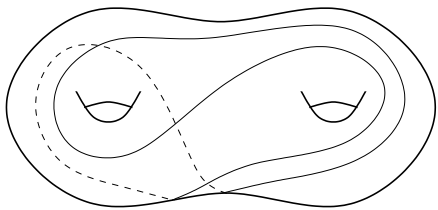}
\caption{The loop $C_0$.}
\label{fig:hyperbolic2}
\end{figure}

Now, we will construct the extension of the $G$-action on $V$ to a
hyperbolic $3$-manifold. As before, let $n$ be any positive integer
divisible by the orders of all the elements of $G$. Assume for the time
being that $g=2$.  We take the same curves $C_1$ and $C_2$ as in the
Seifert-fibered construction, but for $C$ we take the image of the loop
$C_0$ shown in figure~\ref{fig:hyperbolic2} under the homeomorphism of $W$
which is a right-hand twist in a meridinal disc in each of the two
$1$-handles.

\begin{figure}
\includegraphics[width=45 ex]{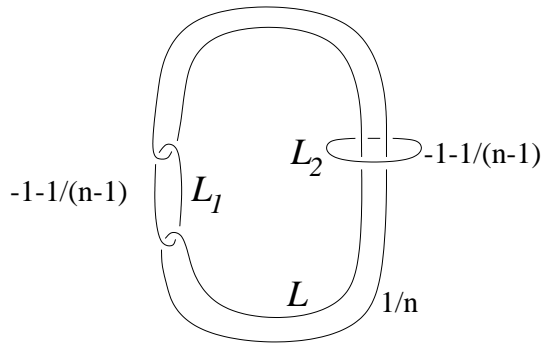}
\caption{Surgery description of $Y$ in the hyperbolic case.}
\label{fig:Whitehead}
\end{figure}

As before, we take the images of the $C_i$ under the $n^{th}$ power of a
Dehn twist about $C$ as the attaching curves, form the closed manifold $Y$,
and use the induced homomorphism $\psi$ to obtain the original $G$-action
as an invariant Heegaard handlebody of a free $G$-action on a covering
space $M$ of $Y$. Let $\ell_1$, $\ell_2$, $L_1$, $L_2$ and $L$ be as
before. We choose the $L_i$ to lie closer to the boundary of $W$ than
$L$. Again, change the attaching curves first by the inverse of the $n$
Dehn twists about $C$, introducing a surgery coefficient of $1/n$ on $L$,
then by the $n-1$ twists about the $\ell_i$, introducing surgery
coefficients $-1/(n-1)$ on the $L_i$. Applying left-hand twists in the
meridian discs of the two $1$-handles of $W$, we move the attaching curves
to the $\ell_i$, obtaining the surgery description of $Y$ shown in
figure~\ref{fig:Whitehead}. This time, the coefficient of $L$ is still
$1/n$, because $L$ has algebraic intersection $0$ with each of the meridian
discs of the handles of $Y$ where the left-hand twists were performed.
The complement of this link is a $2$-fold covering of the complement of the
Whitehead link, so is hyperbolic. By \cite{Thurston}, Dehn surgery on the
link produces a hyperbolic $3$-manifold, provided that one avoids
finitely many choices for the coefficients of each component. So all but
finitely many choices for $n$ yield a hyperbolic $3$-manifold for $Y$, and
hence for $M$.

\begin{figure}
\includegraphics[width=\textwidth]{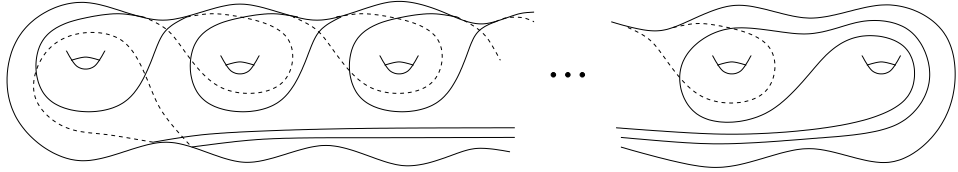}
\caption{The loops $C_0$ for the general hyperbolic construction.}
\label{fig:hyperbolic}
\end{figure}

Finally, to adapt the construction to arbitrary genus, one simply adds more
components to $C_0$ to obtain the $g-1$ circles shown in
figure~\ref{fig:hyperbolic}. In the surgery description for $Y$, the chain
$L\cup L_1$ in figure~\ref{fig:Whitehead} is replaced by a chain of length
$2g-2$, in which $L_1,\ldots\,$, $L_{g-1}$ alternate with components from
$C_0$, and the component $L_g$ links the chain as did $L_2$ in
figure~\ref{fig:Whitehead}.  The $L_i$ have surgery coefficients
$-1-1/(n-1)$ as before, and the components coming from $C_0$ have
coefficient $1/n$, since each had algebraic intersection $0$ with the union
of the meridian discs. The link complement is the $(2g-2)$-fold covering of
the Whitehead link complement, so is hyperbolic, and the argument is
completed as before.

\bibliographystyle{amsplain}

\end{document}